\documentclass[a4paper,12pt]{article}
\usepackage{amsmath,amsthm,amssymb,fontenc,color,textcomp,amsfonts,graphicx,graphics,tabulary}

\begin{document}

\begin{center}{\bf SIMPLIFICATION OF A SYSTEM OF GEODESIC EQUATIONS BY REFERENCE TO CONSERVATION LAWS}\end{center}
\vskip 1.0truecm

\begin{center}{ Uchechukwu Opara }\end{center}
\begin{center}{ Veritas University (VUNA), Abuja, Nigeria \\ e-mail: ucmiop@yahoo.com . }\end{center}
\vskip 1.0truecm
\textbf{Abstract.} \ \ This paper is purposed to exploit prevalent premises for determining analytical solutions to differential equations formulated from the calculus of variations.  We realize this premises from the statement of Emmy Noether's theorem; that every system in which a conservation law is observed also admits a symmetry of invariance (\cite{aj} Page 242, \cite{ab} Pages 60-62).  As an illustration, the infinitesimal symmetries for Ordinary Differential Equations (O.D.E's) of geodesics of the glome are explicitly computed and engaged following identification of a relevant conservation law in action.  Further prospects for analysis of this concept over the same manifold are then presented summarily in conclusion. \\ \\
\textbf{Keywords - } Glome, Geodesics, Calculus of Variations, Conservation Laws, Systems of Differential Equations, Infinitesimal Symmetries. \\ \\

\section {\bf Introduction.} \ \ The 3-sphere, otherwise termed the glome, is a Riemannian manifold at the center of several revolutionary conjectures and advancements in modern mathematical theory.  As prominent examples, one may consider the famous Poincar\'{e} conjecture and the Ricci Flow theorem of Hamilton on closed 3-dimensional manifolds with everywhere positive scalar curvature (\cite{ae} page 128).  The former example challenges an interested mind on meticulous details of differential topology, while the latter is a relatively modern sprout of Pseudo-Riemannian geometry requiring assorted topological and analytical tools.  It may be argued that Hamilton's theorem is conveniently suited towards verification of the Poincar\'{e} conjecture from an assortment of partial vantage points.\\

In this study however, the focus is almost entirely on part of the computational wealth available to be harnessed from conservation laws in action during the course of analogous physical processes on the glome.  Since the choice manifold is of intermediate dimension, the computational work included may be readily confirmed manually, or it could perhaps motivate an exploration of digital software for similar problems.  An advantage of the choice of manifold here is that for many coordinate systems used to capture it, the associated formulated equations have solutions in terms of elementary functions. Moreover, the practical scientific essence of the choice manifold cannot be discarded, as several formulations from the Ricci Flow are formally analogous to heat flow along manifolds and some formulations are adapted to the analysis required in modern cosmological mechanics (\cite{ae} Pages 95, 107).\\

In order to come to terms with a conservation law in play in the course of traversing geodesics of a smooth manifold by a particle without slipping, one may consult the weak formulation process from the stage of differentiation in Banach path spaces. Where the geodesic curve computed in the domain of a relevant coordinate system is parametrized by $t$, we may recall that all such curves satisfy $$\Lambda_{u_i} - \dfrac{d}{dt}(\Lambda_{\overset{\bf .}{u_i}}) = 0,$$ given the arclength functional $$s_{I,t} = \int _I \Lambda (t, u_1 (t), u_2 (t), \cdots , u_{n-1}(t), \overset{\bf .}{u_1}, \overset{\bf .}{u_2}, \cdots, \overset{\bf .}{u}_{n-1})dt$$ for an analytical hypersurface embedded in $\mathbb{R}^n$ and spatial parameters $u_i : 1 \leq i \leq n-1$ (\cite{ad} Page 33). Conservation of the quantity $\left[ \Lambda_{u_i} - \dfrac{d}{dt}(\Lambda_{\overset{\bf .}{u_i}}) \right]$ mentioned above seems to bear analogy to a phenomenon in classical mechanics independent of gravitational influence. 

\section{\bf Method of Formulation} The procedure for computing infinitesimal symmetries accommodated by the system of Euler-Lagrange equations is detailed rather succinctly in the following two theorems.\\ \\
{\underline{\bf Theorem 1}} (\cite{aj} Page 253) A connected group of transformations $G$ acting on $M$ is a variational symmetry group of the functional $s_{I,t}$ if and only if $$ pr^{(j)}{\bf v}(\Lambda) + \Lambda \dfrac{d \xi}{dt} = 0 \ ,$$ for all $(t, u^{(j)}) \in M$ and every infinitesimal generator $$ {\bf v} = \xi (t, u) \frac{\partial}{\partial t} + \sum_{\alpha = 1}^{q}\phi_\alpha (t, u) \dfrac{\partial}{\partial u^\alpha}$$ of $G$.  \\ \\
{\underline{\bf Theorem 2}} (\cite{aj} Page 255) If $G$ is a variational symmetry group of the functional $s_{I,t} = \int _I \Lambda (t, u)dt$, then $G$ is a symmetry group of the associated Euler-Lagrange equations.\\

Consider the Lagrangian computed in \cite{ac} for discovering equations of geodesics for the glome:
$$ s_{I, u_1} = \int _I \sqrt{1 + \cos ^2 u_1 \left( \dfrac{du_2}{du_1} \right) ^2 + \cos ^2 u_1 \cos ^2 u_2 \left( \dfrac{du_3}{du_1} \right) ^2} du_1 := \int _I \Lambda du_1$$
for an appropriate interval $I \subseteq [\frac{-\pi}{2}, \frac{\pi}{2}]$. We hereby reckon with the hyperspherical co-ordinate system,
$$\begin{array}{rlcl}
f : & [\frac{-\pi}{2}, \frac{\pi}{2}]^2 \times [0, 2\pi] & \longrightarrow & S^3 \subset \mathbb{R}^4 \\
{}& (u_1 , u_2 , u_3) & \longmapsto & f(u_1 , u_2 , u_3) = (x_1 , x_2 , x_3 , x_4)
\end{array}$$
$f (u_1 , u_2 , u_3) = (\cos u_1 \cos u_2 \cos u_3 , \cos u_1 \cos u_2 \sin u_3 , \cos u_1 \sin u_2 , \sin u_1)$ .\\ 

For this study, we shall identify with the following renaming of variables: $$(s_{I, u_1}, u_1, u_2, u_3) := (\mathcal{L}, x, y, v).$$  Assuming that $\mathcal{L}$ accommodates an infinitesimal symmetry $${\bf v} = \xi \dfrac{\partial}{\partial x} + \phi \dfrac{\partial}{\partial y} + \eta \dfrac{\partial}{\partial v} \ ,$$ we must equivalently have the following from Theorem 2 above: $$ pr^{(1)}{\bf v}(\Lambda) + \Lambda \dfrac{d \xi}{dx} = 0.$$ We have taken $j=1$ in the statement of the theorem, because no derivative in $\Lambda$ is higher than the first. Now, we reckon that  $pr^{(1)}{\bf v} = {\bf v} + \phi^x \dfrac{\partial}{\partial y_x} + \eta^x\dfrac{\partial}{\partial v_x}$, \  where \\ $\begin{array}{rcl} \phi ^x & = & \mathcal{D}_x(\phi - \xi y_x) + \xi y_{xx} \\ {}&=& \phi_x + \phi_y y_x + \phi_v v_x - [\xi _x y_x + \xi _y (y_x)^2 + \xi_v v_x y_{x}], \end{array}$  \\
$\begin{array}{rcl} \eta ^x & = & \mathcal{D}_x(\eta - \xi v_x) + \xi v_{xx} \\ {}&=& \eta_x + \eta_y y_x +\eta_v v_x - [\xi _x v_x + \xi _y y_x v_x + \xi_v (v_x)^2 ]. \end{array}$ \\ 

$ pr^{(1)}{\bf v}(\Lambda)  = \left( \begin{array}{c} \xi \dfrac{\partial}{\partial x} + \phi \dfrac{\partial}{\partial y} + \eta \dfrac{\partial}{\partial v}+ \\ \phi_x + \phi_y y_x + \phi_v v_x - [\xi _x y_x + \xi _y (y_x)^2 + \xi_v v_x y_{x}]+ \\ \eta_x + \eta_y y_x +\eta_v v_x - [\xi _x v_x + \xi _y y_x v_x + \xi_v (v_x)^2] \end{array} \right) (\Lambda)$ \\ \\
$= \dfrac{\xi}{2\Lambda} (-2(\cos x \sin x) y_x^2 - 2(\cos x \sin x \cos^2 y) v_x^2) \ + \ \dfrac{\phi}{2\Lambda}(-2(\cos^2 x \cos y \sin y) v_x^2)$ \\ $+\dfrac{1}{2\Lambda}(\phi _x + \phi_y y_x + \phi_v v_x - [\xi _x + \xi _y y_x + \xi _v v_x]y_x) (\cos^2 x) . 2y_x$ \\
$+\dfrac{1}{2\Lambda}(\eta _x + \eta_y y_x + \eta_v v_x - [\xi _x + \xi _y y_x + \xi _v v_x]v_x) (\cos^2 x \cos^2 y). 2v_x \ .$

\section{\bf Computational Results}
For fulfillment of admittance of the infinitesimal symmetry criterion derived above from Theorem 2, we must solve the equation: $$\Lambda .pr^{(1)}{\bf v}(\Lambda) + \Lambda^2 .(\xi _x + \xi _y y_x + \xi _v v_x) = 0.$$  We hence evaluate the coefficients of the various uneliminated monomials involved in this equation to zero, as obtained in the table below.
$$\begin{array}{ccl} \underline{MONOMIAL} & \underline{COEFFICIENT} & {} \\ 1 & \xi _x = 0 & (a) \\ y_x & \phi_x \cos^2 x + \xi _ y = 0 & (b) \\ v_x & \eta_x \cos^2 y \cos^2 x + \xi _v = 0 &(c) \\ (y_x)^2 & -\xi \cos x \sin x + \phi _y \cos^2 x = 0 & (d) \\ (v_x)^2 & -\xi \sin x \cos y - \phi \cos x \sin y + \eta_v \cos x \cos y = 0 & (e) \\  y_x v_x & \eta _y \cos^2 x \cos^2 y + \phi _v \cos^2 x = 0 & (f)  \end{array}$$  

We determine the following system of equations from the constraints obtained above in the table of monomial coefficients. $$(a.) \ \xi = \xi (y, v), \ (b.) \ \xi_y = -\phi_x \cos^2 x, \  (c.) \ \xi _v = -\eta_x \cos^2 x \cos ^2 y, \ ( d.) \ \xi = \phi_y \cot x ,$$  
$$ (e.) \ -\xi\sin x \cos y - \phi \cos x \sin y + \eta _v \cos x \cos y = 0 , \ (f.) \ \eta_y \cos^2 y + \phi_v = 0.$$

From $(a.),$ we have $\xi$ to be a function of $y$ and $v$ alone.\\

From $(b.),$ we have $\xi _{yy} = -\phi _{xy}\cos ^2x$, and from (d.),  $$\phi_{xy}\cos ^2x = \xi$$ by differentiating partially with respect to $x$.  Upon comparison with $(b.), \\ \xi _{yy} = \xi$. Hence, we determine that $\xi = \alpha(v)\cos y + \beta(v)\sin y \ .$\\

Substituting this back in $(d.), \ \ \ \phi = \tan x[\alpha(v) \sin y - \beta (v) \cos y] + \gamma (v,x) .$\\ \\
From $(c.),$ we have $$\alpha '(v)\cos y + \beta '(v) \sin y = -\eta _x \cos^2 x \cos ^2y \ \Leftrightarrow$$ $$\eta = -\alpha '(v)\sec y \tan x - \beta '(v)\tan y \sec y \tan x + \delta (y,v) \ .$$
Upon substitution of the above obtained expressions for $\phi$ and $\eta$ in $(f.),$ we realize that $$-2\beta '(v) \sec y \tan x + \delta _y \cos^2 y + \gamma _v = 0$$ after simplification.
Substituting the expressions obtained for $\xi , \eta$ and $\phi$ in $(e.)$, we realize that $$(-\alpha (v) -\alpha ''(v) - \beta ''(v)\tan y)\tan x - \gamma \sin y + \delta _v \cos y = 0 .$$  We shall take $\{k_j\}_{j \in \mathbb{N}}$ to be arbitrary constants in what ensues.  Feasible deductions from the above equation are the following: \\ \\
(I.) $\alpha (v) + \alpha ''(v) = 0 \ \ \Leftrightarrow \ \ \alpha (v) = k_1 \cos v + k_2 \sin v$ \\
(II.) $\beta (v) = k_3$ \\
(III.) $\delta _v = \gamma \tan y$ \\
(IV.) $\gamma = \gamma (v)$ \\ \\
As such, from (III.) and (IV.) we have \  $\delta _y = \kappa (v) \sec ^2y$, where $\kappa = \int \gamma (v) dv$. Substituting these values in the modified equation for $(f.)$, we realize 
$$\int \gamma (v) dv + \gamma _v = 0 \ \Longrightarrow \ \gamma(v) + \gamma _{vv} = 0 \ \Longrightarrow \gamma(v) = k_4 \cos v + k_5 \sin v .$$
The final variable needed to be determined is $\delta:$ $$\delta = \int \gamma dv \tan y \ = \ [k_4 \sin v - k_5 \cos v]\tan y.$$

Upon determining solutions to this system, we observe a general accommodated infinitesimal symmetry of the system of O.D.E's for the glome's geodesics via the hyperspherical coordinate system:\\ \\
${\bf v} = [(k_1 \cos v + k_2 \sin v)\cos y + k_3 \sin y]\dfrac{\partial}{\partial x}$\\ 
$+ [(k_1\cos v + k_2 \sin v)\tan x \sin y - k_3 \tan x \cos y + k_4 \cos v + k_5 \sin v]\dfrac{\partial}{\partial y}$ \\ 
$+ [(k_1\sin v - k_2\cos v)\tan x \sec y + (k_4\sin v - k_5\cos v)\tan y]\dfrac{\partial}{\partial v}.$\\

We can then separate this general symmetry into one-parameter symmetries by the five constants $\{k_i\}_{i=1}^5$ as follows - $$\sum_{i=1}^5 k_i{\chi}_i .$$  By computing the Lie brackets of these accommodated single parameter subgroups, one can then determine stability of the infinitesimal symmetry system, and determine the subgroups of the overall admitted invariance symmetry group of the Lagrangian in question. A sixth accommodated infinitesimal symmetry is actually revealed in this process, being the trivial translation one-parameter group $\left[ \chi _6 = \dfrac{\partial}{\partial v}\right]$. Computation of Lie brackets is particularly instrumental in reducing the associated Euler-Lagrange ordinary differential equation system, by exposing the right \textit{invariants} to be employed together. The element in the $i$'th row and $j$'th column of the table of Lie brackets below is the vector field $[\chi _i , \chi _j]$.\\

Note that the characterization of Lie brackets $$[X, Y] f = X(Y(f)) - Y(X(f))$$ for $X = X^i \frac{\partial}{\partial x_i} , Y = Y^j\frac{\partial}{\partial x_j}$ and any $C^\infty$ function $f$ gives us the formula  
$$ [X, Y] = \sum_i \sum_j \left\{ X^j \frac{\partial Y^i}{\partial x_j} - Y^j \frac{\partial X^i}{\partial x_j} \right\} \frac{\partial}{\partial x_i} .$$

\begin{center}
\begin{tabulary}{ 1.17\textwidth}{ | C |  C   C   C   C   C  C  | }
 \hline
 {} & ${\chi}_1$ & ${\chi}_2$ & ${\chi}_3$ & ${\chi}_4$ & ${\chi}_5$ & ${\chi}_6$ \\
 \hline
 $\chi_1$ & 0 & $-\chi_6$ & $-\chi_4$ & ${\chi}_3$ & 0 & $\chi_2$ \\
${\chi}_2$ & $\chi_6$ & 0 & $-\chi_5$ & 0 & ${\chi}_3$ & $-{\chi}_1$ \\
${\chi}_3$ & $\chi_4$ & $\chi_5$ & 0 & $-\chi_1$ & $-\chi_2$ & 0\\
${\chi}_4$ & $-{\chi}_3$ & 0 & $\chi_1$ & 0 & $ - \chi_6$ & ${\chi}_5$\\
${\chi}_5$ & $0$ & $-{\chi}_3$ & $\chi_2$ & $\chi_6$ & 0 & $-\chi_4$\\
${\chi}_6$ & $-\chi_2$ & ${\chi}_1$ & 0 & $-{\chi}_5$ & $\chi_4$ & 0 \\
\hline 
\end{tabulary}
\end{center} 
\vspace{0.75cm}

We hereby ascertain stability of the accommodated Lie group of infinitesimal symmetries for the formulated geodesic variational problem, with four subgroups listed below. \\
(a.) $\{ \chi _1 , \chi _2 , \chi _ 6\}$\\  
(b.) $\{ \chi _1 , \chi _3 , \chi _ 4\}$\\
(c.) $\{ \chi _4 , \chi _5 , \chi _ 6\}$\\  
(d.) $\{ \chi _2 , \chi _3 , \chi _ 5\}$\\

A significant prospective benefit of the computation done here is an avenue for utilisation of \textit{joint invariants} of any subgroup of the differential equation system in an effort to simplify the simultaneous pair of O.D.E's determined via the Lagrangian. However, in this body of work, we shall make direct reference to the collapsed equation in \cite{ac}:\\  \\ $\overset{ \bf .}{y}\sin x \cos y (k - 2\cos ^2 x \cos ^2 y ) + \overset{ \bf ..}{y}\cos x \cos y ( \cos ^2 x \cos ^2 y - k) +  k \sec x \sin y + k (\overset{\bf .}{y})^2 \cos x \sin y - (\overset{\bf .}{y})^3 \cos ^4 x \sin x \cos ^3 y = 0 \ \ \ \cdots (E = 0) $ \\  \\ derived from substituting the outcome of one Euler-Lagrange equation: $$\Lambda_{v} - \dfrac{d}{dx}(\Lambda_{\overset{\bf .}{v}}) = 0,$$ in the formulation of the second equation in the pair: $$\Lambda_{y} - \dfrac{d}{dx}(\Lambda_{\overset{\bf .}{y}}) = 0.$$  

For emphasis, the constant $k$ in the O.D.E above is in the interval $[0, 1]$ and an overset dot connotes differentiation with respect to the independent variable $x$ of the system.  The reader may confirm that $(E=0)$ accommodates the one-parameter subgroup $[\chi_3 = \sin y \frac{\partial}{\partial x} - \cos y\tan x \frac{\partial}{\partial y}]$  computed earlier, since $[pr^{(2)}\chi _3 (E) = 0]$ whenever $[E = 0]$.  Reckoning with the common benefit of knowledge of an accommodated one-parameter symmetry of an O.D.E, we have foresight that the equation [E = 0] of the second order would be reduced to a significantly simpler equation of lower order, upon replacement of the variables $x$ and $y$ with the canonical coordinates of the symmetry $\chi _3$. Let the canonical coordinates of this symmetry be $(\tau , \omega)$,  whereby $\omega$ is the invariant and $\tau$ the other canonical coordinate.  To determine $\omega$, we solve the equation $$\frac{dx}{siny} = \frac{dy}{-\tan x \cos y} \ ,$$ revealing that $$\cos x . \cos y = [constant] \ .$$  This means we may take $\omega$ to be $\cos x \cos y$, or any smooth function of this term. Moreover, $$\tau = \int \frac{dx}{\sin y}$$ whereby $y$ is expressed in terms of $x$ and $\omega$ (using $\omega = \cos x \cos y$), and $\omega$ is momentarily taken as a constant when evaluating this integral (\cite{ah} Page 24). Thus, we obtain $\tau (x, y) = \arctan ( \cot x \sin y)$ as the second canonical coordinate.   \\

Besides the infinitesimal form of a one-parameter symmetry, the global form also reveals other details present.  The global form $(X(x, y, \lambda ), Y(x, y, \lambda ))$ is determined by integrating the autonomous O.D.E system: $$\frac{dX}{d \lambda} = \sin Y \ , \ \ \frac{dY}{d \lambda} = -\tan X \cos Y$$  subject to the initial constraints $(X, Y)|_{\lambda = 0} = (x, y)$, whereby $x$ and $y$ are temporarily taken as constants in the course of integration.  The solutions obtained here are $$X = \arcsin [ \sin x \cos \lambda - \cos x \sin y \sin \lambda ] \  \ \ \mbox{and} \ \ \ Y = \arctan [ \tan y \cos \lambda + \tan x \sec y \sin \lambda ].$$  

These functions satisfy the customary requirements: $\omega (X, Y) = \omega (x, y)$ \ and \ $\tau (X, Y) = \tau (x, y) + \lambda$.

\section{\bf Discussion}
The equation $[E = 0]$ in focus becomes significantly simpler upon replacement of the initial variables $(x, y)$ with the canonical coordinates $(\tau, \omega)$ of the accommodated symmetry we have engaged. Specifically, this equation impressively reduces to the first order representation: $$\omega'(\tau) = (1-\omega^2) \tan\left[\pm\arccos\sqrt{\frac{\alpha-(\omega^2\cos^2\tau +\sin^2\tau)}{(\frac{k}{\omega^2}-1)(\omega^2\cos^2\tau+\sin^2\tau)}} + \arctan\left(\frac{\omega}{\tan \tau}\right)\right]\ .$$ Another constant $\alpha$ of integration in the first order O.D.E above comes about at a stage of the simplification process from the prior form $[E=0]$.\\

A comparison of this result to what holds for the geodesic equation of the simpler 2-sphere $(S^2)$ reveals some interesting connections. This sphere is obtained as the intersection of the glome $S^3$ with the hyperplane $[v \equiv 0]$. Capturing this manifold with the same spherical coordinate system used above for $S^3$ yields the geodesic differential equation: $$\frac{d^2y}{dx^2} = 2\frac{dy}{dx}\tan x + \left(\frac{dy}{dx}\right)^3\sin x \cos x \ .$$ This equation also accommodates the symmetry $[\chi_3 = \sin y \frac{\partial}{\partial x} - \cos y\tan x \frac{\partial}{\partial y}], $ which was used to simplify the geodesic equation for the glome.  It may be argued that this inherited symmetry of the geodesic equations is based on the fact that $S^2$ is a \textit{totally geodesic} hypersurface of $S^3$.  Simply put, all geodesics of the submanifold $S^2$ are also geodesics of the glome, which apparently corroborates the shared symmetry of their geodesic differential equations.  Determination of totally geodesic hypersurfaces invariably simplifies any theoretical and computational considerations from closed Riemannian manifolds. \\

In general, the geodesic equation for $(S^2)$ admits the 3-dimensional infinitesimal rotation group on Euclidean 3-space (\cite{aa} Page 78), which could contribute otherwise to reduction via invariance transformations for the equation $[E=0]$ we have engaged.  Because none of the admitted symmetries was used in collapsing the required pair of Euler-Lagrange equations to obtain $[E=0]$, the resulting equation may retain invariance under those infinitesimal symmetries in terms of the variables $x$ and $y$. As a major part of what has been established in the previous section, the geodesic equations for $(S^3)$ admit the 6-dimensional infinitesimal rotation group $(SO_4)$ on Euclidean 4-space. Hence, there apparently remains a number of significant consequential properties to be harnessed from this identified symmetry group classification, considering also their explicit functional expressions herein made available.

\section{\bf Conclusion}
Geodesics constitute a class of minimal submanifolds of closed 3-manifolds.  Another class of minimal submanifolds is that of minimal hypersurfaces.  This latter class of 2-dimensional minimal submanifolds are particularly well suited to the study of intrinsic properties of their ambient manifolds in static or evolutionary states (\cite{af} Pages 18 - 27).  As a matter of definite interest, minimal hypersurfaces have been explored severally in academic archives to elucidate details of manifolds' evolution in course of the Ricci Flow.  In the weak formulation of Partial Differential Equations (P.D.E's) for stable minimal hypersurfaces, it is relevant to attempt identification of conservation laws in action.  In this way, Noether's theorem becomes a guarantor of simplification of the associated P.D.E's by accommodated symmetries.\\

Besides probing the inner geometry of 3-manifolds using minimal hypersurfaces, the solitons of Ricci Flow are characterized as equilibrium states of the metrics in the diffusion-reaction equations in process (\cite{ae} Page 4).  This again gives rise to associated conservation laws and an invaluable avenue to engage Noether's theorem.  Since its discovery, this theorem has been at the hub of simplification achievements for systems of differential equations, but the chosen vantage point of the glome has specifically been seen to conceal a plethora of theoretical wealth in this regard.

 \end{document}